\documentstyle[amscd,amssymb,verbatim,12pt]{amsart}
\newcommand{\Area}{\mbox{\rm Area}}

\newcommand{\cone}{\mbox{\rm cone}}
\newcommand{\const}{\mbox{\rm const.}}

\newcommand{\Id}{\mbox{\rm Id}}

\newcommand{\Isom}{\mbox{\rm Isom}}

\newcommand{\pt}{\mbox{\rm pt}}

\newcommand{\R}{{\Bbb R}}

\newcommand{\vol}{\mbox{\rm vol}}
\newcommand{\Z}{{\Bbb Z}}
\theoremstyle{plain}
\newtheorem{definition}{Definition}

\newtheorem{lemma}{Lemma}
\newtheorem{sublemma}{Sublemma}
\newtheorem{theorem}{Theorem}

\numberwithin{equation}{section}
\renewcommand{\rm}{\normalshape}
\begin{document}
\title{Manifolds with Quadratic Curvature Decay and Fast Volume Growth}
\author{John Lott}
\address{Department of Mathematics\\
University of Michigan\\
Ann Arbor, MI  48109-1109\\
USA}
\email{lott@@math.lsa.umich.edu}
\thanks{Research partially supported by NSF grant DMS-0072154 and MSRI}
\date{Feb. 27, 2001}
\maketitle
\begin{abstract}
We give sufficient conditions for a noncompact Riemannian manifold,
which has quadratic curvature decay, to have finite
topological type with ends that are cones over spherical space forms.
\end{abstract}

\section{Introduction} \label{Introduction}
Let $M$ be a complete connected noncompact 
Riemannian manifold with a basepoint $*$.
A natural condition to put on $M$ is that of quadratic curvature decay.
To state this condition,
for $m \in M$ and $r > 0$, let $B_r(m)$ denote the
open distance ball around $m$ of radius $r$ and let $S_r(m) = \partial
\overline{B_r(m)}$ denote the distance sphere around $m$ of radius $r$. 
If $P$ is a $2$-plane in $T_mM$, let $K(P)$ denote the
sectional curvature of $P$. Then $M$ has quadratic curvature decay
if for some $C \: > \: 0$,
\begin{equation} \label{eq1.1}
\limsup_{r \rightarrow \infty} \sup_{m \in S_r(\star), \: P \subset T_mM} 
r^2 \: |K(P)| \: \le \: C.
\end{equation}
Note that (\ref{eq1.1}) is scale-invariant,
in that it is unchanged under a constant rescaling of the Riemannian
metric.

In itself (\ref{eq1.1}) 
does not impose any topological restrictions on $M$, as
any smooth connected manifold admits a complete Riemannian metric satisfying
(\ref{eq1.1}) for some $C$ 
\cite[p. 96]{Gromov (1982)}, \cite[Lemma 2.1]{Lott-Shen (2000)}.
However, with additional assumptions one can obtain restrictions on $M$.
For example, if 
\begin{equation} \label{eq1.2}
\limsup_{r \rightarrow \infty} \sup_{m \in S_r(\star), \: P \subset T_mM} 
r^{2(1+\epsilon)} \: |K(P)| \: < \: \infty
\end{equation} 
for some $\epsilon \: > \: 0$ then
Abresch showed that $M$ has finite topological type, i.e. is homeomorphic
to the interior of a compact manifold-with-boundary
\cite{Abresch (1985)}. For other results on manifolds with
faster-than-quadratic curvature decay, see \cite{Abresch (1985)},
\cite{Greene-Petersen-Zhu (1994)} and
\cite{Petrunin-Tuschmann (1999)}. 

If $M$ has quadratic curvature decay
and a volume growth which is slower than that of the Euclidean space of
the same dimension then topological restrictions on the ends of
$M$ were obtained in a paper of the author with Zhongmin Shen 
\cite{Lott-Shen (2000)}.
Along these lines, we remark that a manifold with quadratic curvature decay and
Euclidean volume growth can have infinite topological type
\cite[Section 2, Example 3]{Lott-Shen (2000)}. Furthermore, even if we assume
finite topological type,
the interior of any connected compact manifold-with-boundary
has a
complete Riemannian metric with quadratic curvature decay and Euclidean volume
growth \cite[Section 2, Example 1]{Lott-Shen (2000)}. Hence the assumptions
of quadratic curvature decay and Euclidean volume growth do not in themselves
give interesting topological restrictions.

In \cite{Lott-Shen (2000)} the question was raised as to what one can say
if one assumes that the constant $C$ in (\ref{eq1.1}) is small enough.
In this paper we give some answers to this question. First,
we show that if the constant $C$ is small enough, if we have pinched
Euclidean volume growth and if $M$ is noncollapsed at infinity
in a suitable sense then $M$ has finite topological type, with ends that
are cones over spherical space forms. 

\begin{theorem} \label{thm1}
Given $n \in \Z^+$ and $c, c^\prime \in \R^+$, there is a constant
$\epsilon \equiv \epsilon(n, c, c^\prime) \: > \: 0$ so that if $M$ is a 
complete connected $n$-dimensional Riemannian manifold with basepoint
$\star$ which satisfies
\begin{equation} \label{eq1.3}
\liminf_{r \rightarrow \infty}
\inf_{m \in S_r(\star)} r^{-n} \: \vol(B_{r/2}(m)) \: \ge \: c
\: \: \: \: \: \text{(noncollapsing)},
\end{equation}
\begin{equation} \label{eq1.4}
c^\prime \: - \: \epsilon \: \le \: \liminf_{r \rightarrow \infty}
r^{-n} \: \vol(B_r(\star)) \: \le \: 
\limsup_{r \rightarrow \infty}
r^{-n} \: \vol(B_r(\star)) \: \le \: c^\prime  \: + \: \epsilon
\: \: \: \: \: \text{(Eucl. vol. growth)}
\end{equation}
and 
\begin{equation} \label{eq1.6}
\limsup_{r \rightarrow \infty} \sup_{m \in S_r(\star), \: P \subset T_mM} 
r^2 \: |K(P)| \: \le \: \epsilon
\: \: \: \: \: \text{(quadratic curvature decay)} 
\end{equation}
then $M$ has finite topological type with ends that are cones over
spherical space forms.
That is, for large $R$, $M - \overline{B_R(\star)}$ is 
homeomorphic to
$(0, \infty) \times Y$ for some closed 
manifold $Y$ which is a union of spherical space forms. Furthermore,
$Y$ has volume $n \: c^\prime$ and the cone over $Y$ satisfies (\ref{eq1.3}).
In particular,
there is a finite number of topological possibilities for $Y$, with
the number depending on $c$ and $c^\prime$.
\end{theorem}

Next, we show that there is a surface of infinite
topological type which admits noncollapsing
metrics of roughly Euclidean volume growth, and
arbitrarily pinched quadratic curvature decay.
The existence of such metrics was pointed out to me by
Bruce Kleiner.

\begin{theorem} \label{thm2}
Given $\epsilon \: > \: 0$, there is a surface of infinite
topological type, equipped with a complete Riemannian metric,
along with constants $c, c_1^\prime, c_2^\prime \: > \: 0$ such that
\begin{equation} \label{eq1.7}
\liminf_{r \rightarrow \infty}
\inf_{m \in S_r(\star)} r^{-2} \: \vol(B_{r/2}(m)) \: \ge \: c
\: \: \: \: \: \text{(noncollapsing)},
\end{equation}
\begin{equation} \label{eq1.8}
c_1^\prime \: \le \: \liminf_{r \rightarrow \infty}
r^{-2} \: \vol(B_r(\star)) \: \le \: 
\limsup_{r \rightarrow \infty}
r^{-2} \: \vol(B_r(\star)) \: \le \: c_2^\prime
\: \: \: \: \: \text{(Euclidean volume growth)}
\end{equation}
and 
\begin{equation} \label{eq1.9}
\limsup_{r \rightarrow \infty} \sup_{m \in S_r(\star), \: P \subset T_mM} 
r^2 \: |K(P)| \: \le \: \epsilon.
\: \: \: \: \: \text{(quadratic curvature decay)} 
\end{equation}
\end{theorem}

Finally, we give a result in which the pinched Euclidean volume growth of
Theorem \ref{thm1} is
replaced by a large-scale convexity assumption.

\begin{definition}
A complete connected Riemannian manifold $M$ with basepoint $*$ is large-scale
pointed-convex if there is a constant $C^\prime \: > \: 0$ such that\\
1. For any normalized minimizing 
geodesic $\gamma \: : \: [a,b] \rightarrow M$ and
any $t \in [0,1]$,
\begin{equation} \label{eq1.10}
d(\gamma(t a \: + \: (1-t) b), \star) \: \le \: t \: d(\gamma(a), \star) 
\: + \: (1-t) \: d(\gamma(b), \star) \: + \: C^\prime  
\end{equation}
and \\
2. For any two normalized minimizing geodesics $\gamma_1, \gamma_2 \: : \:
[0,b] \rightarrow
M$ with $\gamma_1(0) \: = \: \gamma_2(0) \: = \: \star$ and any $t \in [0,1]$,
\begin{equation} \label{eq1.11}
d(\gamma_1(t b), \gamma_2(tb)) \: \le \: t \: d(\gamma_1(b), \gamma_2(b))
\: + \: C^\prime.  
\end{equation}
\end{definition}

Examples of large-scale pointed-convex manifolds are simply-connected
manifolds of nonpositive curvature, and Riemannian manifolds whose
underlying metric spaces are Gromov-hyperbolic 
\cite[Chapitre 2, Pf. of Proposition 25]{Ghys-de la Harpe (1990)}.

\begin{theorem} \label{thm3}
Given $n \: > \: 2$ and $c \in \R^+$, there is a constant
$\epsilon \equiv \epsilon(n, c) \: > \: 0$ with the following property.
Suppose that $M$ is a 
complete connected $n$-dimensional Riemannian manifold with basepoint
which is large-scale pointed-convex and which  satisfies
\begin{equation} \label{eq1.12}
\liminf_{r \rightarrow \infty}
\inf_{m \in S_r(\star)} r^{-n} \: \vol(B_{r/2}(m)) \: \ge \: c
\: \: \: \: \: \text{(noncollapsing)}
\end{equation}
and
\begin{equation} \label{eq1.13}
\limsup_{r \rightarrow \infty} \sup_{m \in S_r(\star), \: P \subset T_mM} 
r^2 \: |K(P)| \: \le \: \epsilon.
\: \: \: \: \: \text{(quadratic curvature decay)} 
\end{equation}
Then $M$ has finite topological type, with ends that are 
cones over spherical space forms.
\end{theorem}

The method of proof of Theorem \ref{thm1} is by contradiction. 
Here is the rough argument. Suppose that
we have a sequence of $n$-dimensional Riemannian manifolds 
$\{M_i\}_{i=1}^\infty$ which 
together provide a counterexample to Theorem \ref{thm1}.
Then each $M_i$ has ``bad'' regions arbitrarily far away from the
basepoint.  By rescaling, we can assume that the unit sphere around the
basepoint in each $M_i$
intersects a bad region. We would like to take a convergent subsequence 
of the $M_i$'s in order to argue by contradiction.  We may not be able to 
take a convergent subsequence in the pointed Gromov-Hausdorff sense, 
as the curvatures may not be uniformly bounded below
at the basepoints.  However, we can always take a pointed ultralimit 
$(X_\omega, \star_\omega)$
(see Section \ref{Ultralimits}). Then any ball in $X_\omega$ away from
the basepoint will be the
Gromov-Hausdorff limit of a subsequence of balls in the 
$M_i$'s.
Under our assumptions, $X_\omega - \star_\omega$ will be $n$-dimensional
and flat with volume growth $V(r) \: = 
c^\prime \: r^n$.
Then $X_\omega$ is a cone over a closed manifold $Y$ which is a union
of spherical space forms. It follows that for
an infinite number of $i$'s, the ``bad'' region in $M_i$ was actually good,
which is a contradiction.

To prove Theorem \ref{thm3} we again form an ultralimit $X_\omega$, which
will have a flat metric on $X_\omega - \star_\omega$ and which
will be pointed-convex. If
${\cal C}$ is a connected component of $X_\omega - \star_\omega$ then its
developing map gives an isometric immersion of the universal cover
$\widetilde{\cal C}$ into $\R^n$.
The convexity is used to show that the developing map is an embedding, with
image $\R^n  - \pt.$, from
which the theorem follows.

The structure of the paper is as follows.  In Section \ref{Ultralimits} 
we recall some
facts about ultralimits of metric spaces.  In Section 
\ref{Proof of Theorem 1} we prove
Theorem \ref{thm1}. In Section 
\ref{Proof of Theorem 2} we prove
Theorem \ref{thm2}. In Section \ref{Proof of Theorem 3} we prove
Theorem \ref{thm3} and make some remarks about its hypotheses.

For background information about Gromov-Hausdorff limits and convergence
results, we refer to \cite{Gromov (1999)} and \cite{Petersen (1997)}.

I thank Bruce Kleiner for discussions and for 
providing some key ideas for this paper.  I
also thank Zhongmin Shen for ongoing discussions, and the referee for
a careful reading and important comments. 

\section{Ultralimits} \label{Ultralimits}

If $\omega$ is a nonprincipal ultrafilter on $\Z^+$ and 
$\{X_i\}_{i=1}^\infty$
is a sequence of metric spaces, let $X_\omega$ be the $\omega$-limit of
the $X_i$'s
(see, for example, 
\cite[Section 3.29]{Gromov (1999)},
\cite[Chapter 9]{Kapovich (2001)} and
\cite[Section 2.4]{Kleiner-Leeb (1997)} for background material).
It is a complete
metric space. An element of
$X_\omega$ has a representative $\{x_i\} \in 
\prod_{i=1}^\infty X_i$. Two such
sequences
$\{x_i\}$ and $\{x_i^\prime\}$ are equivalent
if $\lim_{\omega} d_{X_i}(x_i, x_i^\prime) = 0$. The metric on 
$X_\omega$ is 
\begin{equation} \label{eq2.1}
d_{X_\omega}(\{x_i\}, \{x_i^\prime\}) = 
\lim_\omega d_{X_i}(x_i, x_i^\prime).
\end{equation}

If $\{(X_i, \star_i)\}_{i=1}^\infty$ are pointed metric spaces
then the pointed limit
$(X_\omega, \star_\omega)$ is the subset of $X_\omega$ given by representatives
$\{x_i\}$ such that
$\{d_{X_i}(x_i, \star_i) \}_{i=1}^\infty$ is a bounded sequence.
The basepoint $\star_\omega$
in $X_\omega$ has representative $\{\star_i\}$.
If each $X_i$ is a length space then $X_\omega$ is a length
space and minimizing geodesic segments in $X_\omega$ are ultralimits
of minimizing geodesic segments in $\{X_i\}_{i=1}^\infty$
\cite[Proposition 9.4]{Kapovich (2001)}.

If $X$ is a metric space then 
we let $\cone(X)$ denote the cone on $X$, a pointed
metric space.
 \\ \\
{\bf Example :} Fix $\alpha > 1$.
Take $(X_i, \star_i) = (\R^2,0)$ with Riemannian metric 
\begin{equation} \label{eq2.2}
g_i = i^{-2} \left( 
dr^2 + r^{2 \alpha} d\theta^2 \right)
\end{equation}
on $\R^2 - 0 \cong \R^+ \times \frac{\R}{2 \pi \Z}$. Then by definition,
$(X_\omega, \star_\omega)$ is the asymptotic cone of $X_1$. To describe it,
first, by a change of radial coordinate, $g_i$ is equivalent to
$dr^2 + i^{2\alpha - 2} r^{2\alpha} d\theta^2$. Then by a change of
angular coordinate, $X_i$ consists of 
$\R^+ \times \frac{\R}{i^{\alpha - 1} 2 \pi \Z}$ equipped with the metric
$dr^2 + r^{2\alpha} d\theta^2$,
along with the basepoint $\star_i$.
Put $Y_\omega = \lim_\omega \frac{\R}{i^{\alpha - 1} 2 \pi  \Z}$
(an unpointed limit), which is an 
infinite disjoint union of real lines. (Two points in $Y_\omega$, represented
by sequences $\{y_i\}$ and 
$\{y^\prime_i\}$, lie in the same connected component of $Y_\omega$ if and
only if $\lim_\omega d_{Y_i}(y_i, y^\prime_i) < \infty$.) 
Then $X_\omega$ consists of
$\R^+ \times Y$ with the metric 
$dr^2 + r^{2\alpha} g_{Y_\omega}$, along with the basepoint
$\star_\omega$. The manifolds 
$\{X_i\}_{i=1}^\infty$ have uniform quadratic
curvature decay. Clearly
the sequence $\{(X_i, \star_i)\}_{i=1}^\infty$ is not precompact in the
pointed Gromov-Hausdorff topology.  Nevertheless, in a sense it has
well-defined Gromov-Hausdorff limits away from the basepoint.

For a related  relevant example, take 
$(X_i, \star_i) = (\R^2,0)$ with Riemannian metric 
\begin{equation} \label{eq2.3}
g_i =  
dr^2  +  i^2  r^{2} d\theta^2 
\end{equation}
on $\R^2 - 0 \cong \R^+ \times \frac{\R}{2 \pi \Z}$. 
Put $Y_\omega = \lim_\omega \frac{\R}{i 2 \pi \Z}$.
Then $X_\omega = \cone(Y_\omega)$. There is a flat Riemannian
metric on $X_\omega - \star_\omega$.\\

\section{Proof of Theorem 1} \label{Proof of Theorem 1}

Suppose that the theorem is not true.  Then there is a sequence of
pointed complete connected $n$-dimensional
Riemannian manifolds $\{(M_i, \star_i)\}_{i=1}^\infty$ 
such that :\\
1. Condition (\ref{eq1.3}) is satisfied for each $M_i$.\\ 
2. On $M_i$, we have
\begin{equation} \label{eq3.1}
c^\prime \: - \: \frac{1}{i} \: \le \: \liminf_{r \rightarrow \infty}
r^{-n} \: \vol(B_r(\star_i)) \: \le \: 
\limsup_{r \rightarrow \infty}
r^{-n} \: \vol(B_r(\star_i)) \: \le \: c^\prime
\: + \: \frac{1}{i}.
\end{equation}
3. On $M_i$, we have 
\begin{equation} \label{eq3.2}
\limsup_{r \rightarrow \infty} \sup_{m_i \in S_r(\star_i), \: P_i 
\subset T_{m_i}M_i} 
r^2 \: |K(P_i)| \: \le \: \frac{1}{i}.
\end{equation}
4.a. $M_i$ has infinite topological type or\\
4.b. $M_i$ has an end which has no neighborhood homeomorphic to 
$(0, \infty) \times N$ for any closed manifold $N$ which is a union
of spherical space forms.

Define $\rho_i \in C^0(M_i)$ by $\rho_i(m_i) = d(m_i, \star_i)$.
\begin{lemma} \label{lemma1}
For each $i$, there is a sequence $\{r_{i,j}\}_{j=1}^\infty$ of numbers tending
toward infinity such that
for each $j$, there is a connected component $C_{i,j}$ of
$\overline{B_{4r_{i,j}}(\star_i)} - {B_{r_{i,j}}(\star_i)}$ with the
property that it is
not true that the
map $C_{i,j} \rightarrow
[r_{i,j}, 4r_{i,j}]$, given by restriction of $\rho_i$, defines a topological
fiber bundle whose
fiber is a spherical space form.
\end{lemma}
\begin{pf}
Fix $i$. If the lemma is false then there is a number $R > 0$ so
that for all $r > R$ and for each connected component $C$ of
$\overline{B_{4r}(\star_i)} - {B_{r}(\star_i)}$, the map
$\rho_i \big|_{C} : C \rightarrow [r, 4r]$ defines a 
topological fiber bundle whose fiber
is a spherical space form.
In particular, $C$ is homeomorphic to
$[r, 4r] \times N$ for some spherical space form $N$.

Put $s_1 = R + 1$. Then $\overline{B_{4s_1}(\star_i)} - {B_{s_1}(\star_i)}$
is homeomorphic to $[s_1, 4s_1] \times \coprod_{k \in K} N_k$, where
$K$ is an indexing set and each $N_k$ is a spherical space form.
The restriction of $\rho_i$ to 
$\overline{B_{4s_1}(\star_i)} - {B_{s_1}(\star_i)}$ is given by projection onto
the first factor of $[s_1, 4s_1] \times \coprod_{k \in K} N_k$.
As $\overline{B_{3 s_1}(\star_i)} - 
B_{2 s_1}(\star_i)$ is compact, $K$ must be a finite set.
Let $C_k$ be the connected component of 
$\overline{B_{4s_1}(\star_i)} - {B_{s_1}(\star_i)}$ corresponding to
$[s_1, 4s_1] \times N_k$. Put $s_2 =
3s_1$. There is a connected component $C_k^\prime$ of 
$\overline{B_{4s_2}(\star_i)} - {B_{s_2}(\star_i)}$ which intersects
$C_k$. We know that it is homeomorphic to $[s_2, 4s_2] \times 
N^\prime$ for some spherical space form $N^\prime$, with the
restriction of $\rho_i$ to $C_k^\prime$ given by projection onto the
first factor of $[s_2, 4s_2] \times 
N^\prime$.
Then $N^\prime \: = \: N_k$.
Thus $C_k \cup C_k^\prime$ is homeomorphic to $[s_1, 4s_2] \times N_k$ and
extends $C_k$. As each connected component of  
$\overline{B_{4s_2}(\star_i)} - {B_{s_2}(\star_i)}$ intersects
$\overline{B_{4s_1}(\star_i)} - {B_{s_1}(\star_i)}$, we see that
$\overline{B_{4s_2}(\star_i)} - {B_{s_1}(\star_i)}$ is homeomorphic to
$[s_1, 4s_2] \times \coprod_{k \in K} N_k$.
Taking $s_3 =
3s_2$ and continuing the process, we obtain that
$M_i - \overline{B_{R+1}(\star_i)}$ is homeomorphic to
$(0, \infty) \times \coprod_{k \in K} N_k$.
\end{pf}

With reference to Lemma \ref{lemma1}, (\ref{eq1.3}), (\ref{eq3.1})
and (\ref{eq3.2}), we can find a
sequence $R_i = r_{i,j(i)}$ tending towards infinity such that\\
1. For $r \: > \: \frac{1}{i}$,
\begin{equation} \label{eq3.3}
\inf_{m_i \in S_{R_i r}(\star_i)} (R_i  r)^{-n} \: \vol(
B_{R_i r/2}(m_i)) \: \ge \: c \: - \: \frac{1}{i},
\end{equation}
\begin{equation} \label{eq3.4}
c^\prime \: - \: \frac{2}{i} \: \le \: (R_i r)^{-n} \:
\vol(B_{R_i r}(\star_i)) \: \le \: c^\prime \: + \: \frac{2}{i},
\end{equation}
and
\begin{equation} \label{eq3.5}
\sup_{m_i \in S_{R_i r}(\star_i), \: P_i \subset T_{m_i}M_i} (R_i  r)^{2} \: 
|K(P_i)| \: \le \: \frac{2}{i}.
\end{equation}
\noindent
2. There is a connected component $C_i$ of 
$\overline{B_{4R_i}(\star_i)} - {B_{R_i}(\star_i)} \subset M_i$ with the
property that it is not true that the map $C_i \rightarrow
[R_i, 4R_i]$, given by restriction of $\rho_i$, defines a topological
fiber bundle whose
fiber is a spherical space form.

Let $X_i$ be $M_i$ with the rescaled metric $g_{X_i} = 
(2 R_i)^{-2} g_{M_i}$. Define $\mu_i \in C^0(X_i)$ by
$\mu_i(x_i) = d(x_i, \star_i)$.
Let $(X_\omega, \star_\omega)$ be the 
$\omega$-limit of $\{(X_i, \star_i)\}_{i=1}^\infty$.
\begin{lemma} \label{lemma2}
$X_\omega - \star_\omega$ is
a flat $n$-dimensional manifold.
\end{lemma}
\begin{pf}
Given $x_\omega \in X_\omega - \star_\omega$,
put $D = d(x_\omega, \star_\omega)$.
Then $D > 0$. Choose a representative
$\{x_i\} \in \prod_{i=1}^\infty X_i$ of $x_\omega$. 
For any $\epsilon > 0$, there is a subset $W \subset \Z^+$ of full
$\omega$-measure such that for all $i \in W$,
\begin{equation} \label{eq3.6}
|d(x_i, \star_i) - D| < \epsilon.
\end{equation}
If $i \in W$ put $y_i = x_i$ and if $i \notin W$, choose 
$y_i \in S_D(\star_i) \subset X_i$. Then $\lim_\omega d_{X_i}(x_i, y_i) = 0$
and so $\{y_i\}$ also represents $x_\omega$. Thus in replacing $\{x_i\}$ by
$\{y_i\}$, we may assume that $d(x_i, \star_i) \in (D-\epsilon, D+\epsilon)$
for all $i \in \Z^+$. Take $\epsilon \in  \left( 0, \frac{D}{10} \right)$.

Due to the rescaling used to define $X_i$, 
for all $r \: > \: \frac{1}{i}$,
\begin{equation} \label{eq3.7}
\inf_{x_i \in S_r(\star_i)} r^{-n} \: \vol(B_{r/2}(x_i)) \: \ge \: c \: - \:
\frac{1}{i}
\end{equation}
\begin{equation} \label{eq3.8}
c^\prime \: - \: \frac{2}{i} \: \le \: r^{-n} \:
\vol(B_{r}(\star_i)) \: \le \: c^\prime \: + \: \frac{2}{i},
\end{equation}
and
\begin{equation} \label{eq3.9}
\sup_{x_i \in S_r(\star), \: P_i \subset T_{x_i}X_i} 
r^2 \: | K(P_i) | \: \le \: \frac{2}{i}.
\end{equation}
Equation (\ref{eq3.9}) gives a uniform lower
bound on the sectional curvatures of 
$\{\overline{B_{4D/5}(x_i)}\}_{i=1}^\infty$.  It follows that
the closed balls
$\{\overline{B_{3D/4}(x_i)}\}_{i=1}^\infty$ are precompact in the
pointed Gromov-Hausdorff topology \cite[Theorem 2.2, Fact 4]{Petersen (1997)}.
To be precise,  \cite[Theorem 2.2, Fact 4]{Petersen (1997)} deals with
pointed Gromov-Hausdorff precompactness in the case of complete manifolds.
However, in view of the definition of pointed Gromov-Hausdorff precompactness,
the same argument applies to the distance balls.

\begin{sublemma} \label{sublemma1}
$\overline{B_{3D/4}(x_\omega)}$ is a limit point of
$\{\overline{B_{3D/4}(x_i)}\}_{i=1}^\infty$ in the 
pointed Gromov-Hausdorff topology. 
\end{sublemma}
\begin{pf}
The proof is similar to that of \cite[Lemma 2.4.3]{Kleiner-Leeb (1997)}.
By precompactness,
for any $\delta > 0$ there is a number $J$ such that for each $i$, there
is a $\delta$-net $\{x_{i,j}\}_{j=1}^J$ in 
$\overline{B_{3D/4}(x_i)}$, with $x_{i,1} = x_i$.
Let $x_{\omega, j} \in X_\omega$ be represented
by the sequence $\{x_{i,j}\}$. 
In particular, $x_{\omega, 1} = x_\omega$.
We claim that $\{x_{\omega, j}\}_{j=1}^J$ is
a $\delta$-net in $\overline{B_{3D/4}(x_\omega)}$. First,
\begin{equation} \label{eq3.10}
d_{X_\omega}(x_{\omega,j}, x_\omega) = 
\lim_\omega d_{X_i}(x_{i,j}, x_i) \le 3D/4,
\end{equation}
so $x_{\omega, j} \in
\overline{B_{3D/4}(x_\omega)}$. Next, given $y_\omega = \{y_i\} \in
\overline{B_{3D/4}(x_\omega)}$, 
for $j \in \{1, \ldots, J\}$ put
\begin{equation} \label{eq3.11}
U_j = \{i : d_{X_i}(x_{i,j}, y_i) \le \delta\}.
\end{equation}
As $\Z^+ = \bigcup_{j=1}^J U_j$, there is some $j$ so that $U_j$ has full
$\omega$-measure. Then for this $j$,
$d_{X_\omega}(x_{\omega,j}, y_\omega) = 
\lim_\omega d_{X_i}(x_{i,j}, y_i) \le \delta$. Thus 
$\{x_{\omega, j}\}_{j=1}^J$ is a $\delta$-net in 
$\overline{B_{3D/4}(x_\omega)}$.

From the definition of $d_{X_\omega}$, there is a subset $W \subset \Z^+$ of full
$\omega$-measure such that for all $i \in W$ and $j,k \in \{1, \ldots, J\}$,
\begin{equation} \label{eq3.12}
\big| d_{X_\omega}(x_{\omega,j}, x_{\omega,k}) - 
d_{X_i}(x_{i,j}, x_{i,k}) \big| < \delta.
\end{equation}
For any $i \in W$, it follows as in the proof of
\cite[Proposition 3.5(b)]{Gromov (1999)} that the pointed Gromov-Hausdorff 
distance between $\overline{B_{3D/4}(x_\omega)} \subset X_\omega$ and 
$\overline{B_{3D/4}(x_i)} \subset X_i$ is at
most $2 \delta$. This proves the sublemma.
\end{pf}

From (\ref{eq3.7}), for $i$ sufficiently large,
\begin{equation} \label{eq3.13}
\vol (B_{3D/5}(x_i)) \: \ge \: \vol (B_{d(x_i, \star_i)/2}(x_i)) \: \ge \:
\frac{1}{2} \: c \: d(x_i, \star_i)^n \: \ge \: \frac{1}{2} \:  
c \: (9D/10)^n. 
\end{equation}
Hence we are in the noncollapsing situation and so
from \cite[Corollary 2.3, Lemma 3.4 and Theorem 4.1]{Petersen (1997)}, 
$X_\omega - \star_\omega$ has a flat $n$-dimensional Riemannian metric.
\end{pf}

From Sublemma \ref{sublemma1} 
and \cite[Theorem 2.2]{Petersen (1997)},
there is an infinite subset $S \subset \Z^+$ such that
$\overline{B_{3D/5}(x_\omega)}$ is actually the limit of
$\{\overline{B_{3D/5}(x_i)}\}_{i \in S}$ in the 
$C^{1,\sigma}$-topology for any $\sigma \in (0, 1)$.
Given $\alpha > 1$ and $r > 0$, put $A_i(\alpha r, r) = 
\overline{B_{\alpha r}(\star_i)} - B_{r}(\star_i) \subset X_i$ and
$A_\omega(\alpha r, r) = 
\overline{B_{\alpha r}(\star_\omega)} - B_{r}(\star_\omega) \subset 
X_\omega$. From (\ref{eq3.8}), 
\begin{equation} \label{eq3.14}
\lim_{i \rightarrow \infty} \frac{\vol(A_i(\alpha r, r))}{r^n} =
(\alpha^n - 1) \: c^\prime.
\end{equation}
By abuse of notation, we write
$\vol(B_r(\star_\omega))$ for  $\vol(B_r(\star_\omega) - \star_\omega)$.

\begin{lemma} \label{lemma3}
For all $r \: > \: 0$,
\begin{equation} \label{eq3.15}
\vol(A_\omega(\alpha r, r)) = (\alpha^n - 1) \: c^\prime \: r^n
\end{equation}
and
\begin{equation} \label{eq3.16}
\inf_{x_\omega \in S_r(\star_\omega)} r^{-n} \: 
\vol(B_{r/2}(x_\omega)) \: \ge \: c.
\end{equation}
\end{lemma}
\begin{pf}
Given $x_\omega \in 
A_\omega(\alpha r, r)$, let $\{x_i\}_{i=1}^\infty$ be as in the proof of
Lemma \ref{lemma2}. For
$\epsilon > 0$ sufficiently small, 
the method of proof of Sublemma \ref{sublemma1} shows that 
$\overline{B_{\epsilon}(x_\omega)}$ is the pointed Gromov-Hausdorff limit of a
sequence of $\epsilon$-balls $\{ 
\overline{B_{\epsilon}(x_i)} \}_{i \in S}$.
From the Vitali covering theorem \cite[Theorem 2.8]{Mattila (1995)},
if $\vol(A_\omega(\alpha r, r)) \: < \: \infty$ then 
for any $\delta \: > \: 0$ there is a finite number of disjoint closed
metric balls $\{\overline{B(x_{\omega, j}, r_j)}\}_{j=1}^J$ contained in 
$A_\omega(\alpha r, r)$
such that 
\begin{equation} \label{add1}
\sum_{j=1}^J \vol( \overline{B(x_{\omega, j}, r_j)}) \: \ge \:
\vol(A_\omega(\alpha r, r)) \: - \: \delta,
\end{equation}
while if $\vol(A_\omega(\alpha r, r)) \: = \: \infty$ then for
any $\Delta \: > \: 0$, there is a finite number of disjoint closed
metric balls $\{\overline{B(x_{\omega, j}, r_j)}\}_{j=1}^J$ contained in 
$A_\omega(\alpha r, r)$
such that 
\begin{equation} \label{add11}
\sum_{j=1}^J \vol( \overline{B(x_{\omega, j}, r_j)}) \: \ge \:
\Delta.
\end{equation}
(Note that $A_\omega(\alpha r, r)$ could {\it a priori} have an infinite
number of connected components.)
The
$C^{1,\sigma}$ metric convergence implies that 
for any $\epsilon \: > \: 0$
and for an infinite number of $i$'s, there
are disjoint closed metric balls 
$\{\overline{B(x_{i, j}, r_j)}\}_{j=1}^J$ contained in $A_i(\alpha r, r)$ with
\begin{equation} \label{add2}
\sum_{j=1}^J \vol( \overline{B(x_{\omega, j}, r_j)}) \: \le \:
\sum_{j=1}^J \vol( \overline{B(x_{i, j}, r_j)}) \: + \: \epsilon \: \le \:
\vol(A_i(\alpha r, r)) \: + \: \epsilon.
\end{equation}
Equations (\ref{eq3.14}), (\ref{add11}) and (\ref{add2}) imply that in fact 
$\vol(A_\omega(\alpha r, r)) \: < \: \infty$. Then
equations (\ref{eq3.14}), (\ref{add1}) and (\ref{add2}) imply that
$\vol(A_\omega(\alpha r, r)) \: \le \: (\alpha^n \: - \: 1) \: c^\prime \:
r^n \: + \: \delta \: + \: \epsilon$. As $\delta$ and $\epsilon$ are 
arbitrary, we obtain that
\begin{equation} \label{eq3.17}
\vol(A_\omega(\alpha r, r)) \le (\alpha^n - 1) \: c^\prime \: r^n.
\end{equation}

From (\ref{eq3.7}), the lower curvature bound and the Bishop-Gromov
inequality \cite[Lemma 5.3.bis]{Gromov (1999)}, for large $i$
we obtain a lower
bound on $\vol(B_{\epsilon}(x_i))$ in terms of 
$\epsilon$, $\alpha$, $r$ and $c$. Using the $C^{1,\sigma}$ metric
convergence, we obtain a
lower bound on $\vol(B_{\epsilon}(x_\omega))$ in terms of
$\epsilon$, $\alpha$, $r$ and $c$. We then obtain an upper bound on the number
of elements in a maximal $2 \epsilon$-separated net in 
$A_\omega(\alpha r, r)$.
As the $4\epsilon$-balls with centers at the netpoints cover
$A_\omega(\alpha r, r)$, it follows that $A_\omega(\alpha r, r)$ is
compact.  Then
$A_\omega(\alpha r, r)$ is the Gromov-Hausdorff limit of
a subsequence of $\{A_i(\alpha r, r)\}_{i=1}^\infty$. It follows from
the $C^{1,\sigma}$ metric convergence that 
\begin{equation} \label{eq3.18}
\vol(A_\omega(\alpha r, r)) = 
\lim_{i \rightarrow \infty} {\vol(A_i(\alpha r, r))} = 
(\alpha^n - 1) \: c^\prime \: r^n.
\end{equation} 
Equation (\ref{eq3.16}) follows from (\ref{eq3.7}) and 
the $C^{1,\sigma}$ metric convergence.
\end{pf}

Hence $\vol(B_r(\star_\omega)) \: = \: c^\prime \: r^n$.
As $A_\omega(\alpha r, r)$ is compact, we can now use the analysis 
of manifolds that are flat outside of a compact set, as given in
\cite{Eschenburg-Schroeder (1987)}. For simplicity suppose that
$X_\omega \: - \: \star_\omega$ is connected; the general case is
similar.  Suppose that $n \: > \: 2$.
From \cite{Eschenburg-Schroeder (1987)}, the complement of some bounded
set in $X_\omega$ is isometric to the complement of a bounded set in 
$\R^n/F$, for some finite group
$F \subset O(n)$ that acts freely on $S^{n-1}$. For
$r_0$ large, we identify $S_{r_0}(\star_\omega)$ with a hypersurface in
$\R^n/F$. Then for $r \: < \: r_0$, $S_{r}(\star_\omega)$ is the 
result of (possibly) making identifications on 
the equidistant set with signed distance $r \: - \: r_0$ from 
$S_{r_0}(\star_\omega)$. 
We know that 
\begin{equation} \label{area1}
\Area(S_r(\star_\omega)) \: = \: n \: c^\prime \: r^{n-1}.
\end{equation}
As this is analytic in $r$, 
it follows that there are in fact no identifications made,
and $S_{r_0}(\star_\omega)$ is convex when lifted to $\R^n$.
If $r_0$ is large enough, we may assume that
$S_{r_0}(\star_\omega)$ is $C^1$-smooth with 
measurable principal curvature functions
$\{h_j\}_{j=1}^{n-1}$. 
For $r$ near $r_0$, the tube formula gives
\begin{equation} \label{area2}
n \: c^\prime \: r^{n-1} \: = \: 
\int_{S_{r_0}(\star_\omega)} \prod_{j=1}^{n-1} (1 \: + \:
h_j \: (r \: - \: r_0)) \: d\vol.
\end{equation}
By analyticity, (\ref{area2}) is true for all $r$.
As in the proof of the Bishop-Gromov inequality, for $r \in (0, r_0)$,
\begin{equation} \label{area3}
\Area(S_r(\star_\omega)) \: = \: 
\int_{S_{r_0}(\star_\omega)} \chi_r \: \prod_{j=1}^{n-1} 
(1 \: + \:
h_j \: (r \: - \: r_0)) \: d\vol,
\end{equation}
where $\chi_r$ is the characteristic function of the set of points on
$S_{r_0}(\star_\omega)$ whose normal rays are distance-minimizing down to
$S_{r}(\star_\omega)$. It follows from (\ref{area1}), (\ref{area2}) 
and (\ref{area3})
that for all $r \in (0, r_0)$, $\chi_r \: = \: 1$ and
$1 \: + \: h_j \: (r \: - \: r_0) \: > \: 0$ for all $j$.
Equation (\ref{area2}), for small $r$, now implies that for all $j$, 
$h_j \: = \: \frac{1}{r_0}$.  Then for all $r \: > \: 0$,
$S_r(\star_\omega)$ can be identified with
the sphere of distance $r$ from the vertex of
$\R^n/F$. Hence $X_\omega - \star_\omega$ is
a cone over a spherical space form.  If $n \: = \: 2$ then one
can apply a similar argument, using the results of 
\cite{Eschenburg-Schroeder (1987)} in this case.

Hence $X_\omega$ is a cone over a finite
union $Y$ of spherical space forms
with total volume $n \: c^\prime$. 
Let $C_i \subset X_i$ be as above.  Choose $c_i \in C_i \cap S_1(\star_i)$.
Let $c_\omega \in X_\omega$ be the point represented by $\{c_i\}$.
Consider the connected component ${\cal C}$ of $X_\omega - \star_\omega$
which contains $c_\omega$.
Define $\nu_\omega \in C^0({\cal C})$ by
$\nu_\omega(x_\omega) = d_{X_\omega}(x_\omega, \star_\omega)$.

Consider the closed annulus 
$A = \overline{B_4(\star_\omega)} - B_{\frac14}(\star_\omega)$
in ${\cal C}$. It is compact. Given
$\epsilon \in \left( 0, \frac{1}{100} \right)$, 
choose a finite $\epsilon$-net ${\cal N} = \{a_{\omega,j}\}_{j=1}^J$ in $A$,
with $a_{\omega,1} = c_\omega$.
For each $j$,
choose a sequence $\{a_{i,j}\}$ which represents $a_{\omega,j}$, with
$a_{i,j} \in X_i$ and $a_{i,1} = c_i$. 
As in the proof of Lemma \ref{lemma2}, we may assume that
$a_{i,j} \in \overline{B_{4+\epsilon}(\star_i)} - 
B_{\frac14 - \epsilon}(\star_i)$.
By the definition of $d_{X_\omega}$, 
there is a subset $S_0 \subset \Z^+$ of
full $\omega$-measure such that if $i \in S_0$ then
for all $j \in \{1, \ldots, J\}$,
\begin{equation} \label{eq3.19}
\big| d_{X_i}(a_{i,j}, \star_i) - d_{X_\omega}(a_{\omega,j}, \star_\omega) 
\big| < \epsilon.
\end{equation}

Consider the closed subsets 
$\left\{ \bigcup_{j=1}^J \overline{B_{1/8}(a_{i,j})} \right\}_{i \in S_0}$
of $\{X_i\}_{i \in S_0}$. From (\ref{eq3.9}),
they form a precompact set in the multipointed Gromov-Hausdorff topology,
where the multibasepoint of 
$\bigcup_{j=1}^J \overline{B_{1/8}(a_{i,j})}$
is the ordered set
$\{a_{i,j}\}_{j=1}^J$ and by ``multipointed Gromov-Hausdorff topology''
we mean the analog of the pointed Gromov-Hausdorff topology, in which
all of the maps in the definitions respect the multibasepoints.
Put 
\begin{equation} \label{eq3.20}
F = \bigcup_{j=1}^J \overline{B_{1/8}(a_{\omega,j})} \subset
\overline{B_{10}(\star_\omega)} -
B_{\frac{1}{10}}(\star_\omega) \subset {\cal C}.
\end{equation}
As in the proof of Sublemma \ref{sublemma1}, $F$ is a limit point of
$\left\{
\bigcup_{j=1}^J \overline{B_{1/8}(a_{i,j})} \right\}_{i \in S_0}$ 
in the multipointed 
Gromov-Hausdorff topology. Then there is a subsequence of 
$\left\{
\bigcup_{j=1}^J \overline{B_{1/8}(a_{i,j})} \right\}_{i \in S_0}$ which
converges in the multipointed
$C^{1,\sigma}$-topology to $F$.
In particular,
there is an infinite subset $S_1 \subset S_0$ such that
if $i \in S_1$ then there is a $C^{2,\sigma}$-regular diffeomorphism
$ \pi_i : 
\bigcup_{j=1}^J \overline{B_{1/8}(a_{i,j})} \rightarrow 
F$
with $\pi_i(a_{i,j}) = 
a_{\omega,j}$.

For $i \in S_1$, let $g_i$ denote
the corresponding 
Riemannian metric on $F$, pulled back from $X_i$ via
$\pi_i^{-1}$.
Then
$\{g_i\}_{i \in S_1}$ converges to the $g_\omega \big|_F$ 
in the $C^{1,\sigma}$-topology.
Taking $\epsilon \in  \frac{\Z^+}{100}$ and
doing a diagonal argument, we obtain a sequence parametrized by 
$k \in \Z^+$ of 
\\
1. $\frac{1}{100k}$-nets
${\cal N}_k = \{a_{\omega,j}\}_{j=1}^{J_k}$ in $A$,\\
2. Sets
$F_k = \bigcup_{j=1}^{J_k} \overline{B_{1/8}(a_{\omega,j})}$
and \\
3. $C^{2,\sigma}$-regular diffeomorphisms
$\pi_k : 
\bigcup_{j=1}^{J_k} \overline{B_{1/8}(a_{i_k,j})} \rightarrow F_k$ with 
$\pi_k(a_{i_k,j}) = 
a_{\omega,j}$\\
such that \\
4. $\lim_{k \rightarrow \infty} g_{k} \big|_A = g_\omega \big|_{A}$ 
in the $C^{1,\sigma}$-topology and \\
5.
\begin{equation} \label{eq3.21}
\sup_{y \in {\cal N}_k} \big| \nu_{k}(y) - \nu_\omega(y) \big| <
\frac{1}{k},
\end{equation}
where $\nu_k = \mu_{i_k} \circ \pi_{k}^{-1} \in C^0(F_k)$.
Here 
$g_k$ is the pullback of the Riemannian metric from $X_{i_k}$, using 
$\pi_k^{-1}$.
By the Arzela-Ascoli theorem,
it follows from 4. that there is a subsequence of
$\{ \nu_k \big|_{A} \}_{k=1}^\infty$ which converges in the Lipschitz
topology.  Relabelling this subsequence as 
$\{ \nu_k \big|_{A} \}_{k=1}^\infty$, it follows from (\ref{eq3.21}) that
$\lim_{k \rightarrow \infty} \nu_k \big|_{A}
 = \nu_\omega \big|_{A}$.
For large $k$, we will identify $C_{i_k}$ with the connected
component of
$\nu_k^{-1} \left( \left[ \frac{1}{2}, 2 \right] \right) \subset A$ 
containing $c_\omega$.

Let $r$ be the coordinate on ${\cal C}$ given by the distance from
$\star_\omega$ and let $Z = - \: 
\frac{d}{dr}$ be the corresponding (smooth) vector field
on ${\cal C}$. Clearly $Z$ is transversal to $\nu_\omega$ in the sense of
\cite{Greene-Shiohama (1981)}. Then for large $k$, $Z \big|_{A}$ is 
transversal
to $\nu_k \big|_A$. By flowing along $Z$ from $\nu_k^{-1}(2)$ to 
$\nu_k^{-1} \left( \frac{1}{2} \right)$ and using the arguments of
\cite{Greene-Shiohama (1981)}, it follows that the map $\nu_k : 
\nu_k^{-1} \left( \left[ \frac{1}{2}, 2 \right] \right) \rightarrow 
\left[ \frac{1}{2}, 2 \right]$ defines a topological fiber bundle.
By further flowing along $Z$ down to 
$S_{1/100}(\star_\omega) \subset {\cal C}$,
it follows that the fiber of the bundle is homeomorphic to a connected
component of
$Y$. Then $C_{i_k}$ is the total space of this fiber bundle, which 
contradicts the construction of $\{C_i\}_{i=1}^\infty$.

For given $c$ and $c^\prime$, $Y$ has volume $n \: c^\prime$ and
$\cone(Y)$ satisfies
the lower local volume bound (\ref{eq3.16}).
It follows that there is
an upper bound in terms of $c$ and $c^\prime$ on the number of components 
of $Y$, and a finite number of possible diffeomorphism types for
each component. $\square$

\section{Proof of Theorem 2} \label{Proof of Theorem 2}

The underlying basis for the result is the fact that for $\beta \in (0,1)$,
there is a flat
$2$-dimensional cone surface with one cone point, of total angle
$2 \pi (1 \: + \: \beta)$, and one open end, with cone angle 
$2 \pi (1 \: - \: \beta)$. Because of this fact, it is plausible that one can
construct a sequence of surfaces as in the statement of the theorem with
the property that when one takes an ultralimit as $\epsilon \rightarrow 0$,
one obtains this flat cone surface.

This suggests constructing the surface of the theorem to have a self-similar
structure of the form
\begin{equation} \label{add3}
M \: = \: D^2 \: \cup_{S^1} P \: \cup_{S^1}  (C \cdot P) \cup_{S^1}
\ldots \cup_{S^1} (C^k \cdot P) \cup_{S^1} \ldots
\end{equation}
Here $P$, the basic building block, is the gluing $N_1 \cup_{S_1} N_2$ of
two compact surfaces-with-boundary $N_1$ and $N_2$ along a circle.
The surface $N_1$ will be the above cone surface truncated both near the
cone point and near infinity. Topologically $N_1$
will be a torus with two balls removed, 
equipped with a flat metric. Then the surface $N_2$ will be an annulus that
attaches $N_1$ and a rescaled version $C \cdot N_1$, for an 
appropriate constant $C$.

To write this in detail,
let $T^2$ denote the $2$-torus equipped with an arbitrary but fixed
complex structure with local complex coordinate $z$, 
and flat Riemannian metric $|dz|^2$.  
Let $f$ be a meromorphic function on $T^2$ with one zero, at $p_0 \in
T^2$, and one pole, at $p_\infty \in T^2$. Fix $\beta \in (0,1)$ and
put $g \: = \: |f(z)|^{2 \beta} \: |dz|^2$, a Riemannian metric on
$T^2 \: - \: \{p_0, p_\infty\}$. In general, a metric 
$e^{2 \phi} \: |dz|^2$ has Gaussian curvature
$- \: e^{-2 \phi} \: (\partial_x^2 \: + \: \partial_y^2) \phi$. As 
$\ln |f|$ is harmonic, it follows that $g$ is flat.
As a metric on $T^2$, it has a cone
point at $p_0$ with total angle $2 \pi (1 \: + \: \beta)$ (i.e.
angle excess $2 \pi \beta$) and an open cone near $p_\infty$ with
cone angle $2 \pi (1 \: - \: \beta)$.  The end of 
$T^2 \: - \: \{p_0, p_\infty\}$ approaching 
$p_0$ has a neighborhood $U_0$ with the
metric $ds^2 \: + \: (1 \: + \: \beta)^2 \: s^2 \: d\theta^2$ for
$s \in (0, \delta_0)$, and 
the end of 
$T^2 \: - \: \{p_0, p_\infty\}$ approaching 
$p_\infty$ has a neighborhood $U_\infty$ with the
metric $dt^2 \: + \: (1 \: - \: \beta)^2 \: t^2 \: d\theta^2$ for
$t \in (\delta_\infty, \infty)$. We take $\delta_\infty \: > \: 1$.
Put $N_1 \: = \: (T^2 \: - \: \{p_0, p_\infty\}) \: - \: U_0 \: - \: 
U_\infty$. It is a compact surface-with-boundary whose boundary
circles $\partial_0(N_1)$ and $\partial_\infty(N_1)$ have lengths 
$2 \pi (1 \: + \: \beta) \delta_0$ and
$2 \pi (1 \: - \: \beta) \delta_\infty$, respectively.  
If $C$ is a positive constant, we denote by $C \cdot N_1$ the
Riemannian manifold obtained by rescaling the Riemannian metric on 
$N_1$ by $C^2$, i.e. multiplying the lengths by $C$.

For $\epsilon$ a small positive number, we now wish to construct a metric
\begin{equation} \label{eq4.1}
ds^2 \: = \: dr^2 \: + \: f^2(r) \: d\theta^2
\end{equation}
on an annulus $N_2 \: = \: [0, R] \times S^1$ with
\begin{equation} \label{eq4.2}
f(r) \: = \: c_1 \: (1 \: + \: r)^{- \epsilon} \: + \:
c_2 \: (1 \: + \: r)^{1 \: + \: \epsilon}, 
\end{equation}
so that $N_1$ glues isometrically to $N_2$ to first order,
with $\{0\} \times S^1$ gluing to 
$\partial_\infty(N_1)$, and $N_2$ glues isometrically to $C \cdot N_1$
to first order, with
$\{R\} \times S^1$ gluing  to 
$C \cdot \partial_0(N_1)$, for some $C \: > \: 1$.
These conditions become
\begin{align} \label{eq4.3}
c_1 \: + \: c_2 \: & = \: (1 \: - \: \beta) \:  \delta_\infty, \notag \\
- \: \epsilon \: c_1 \: + \: (1 \: + \: \epsilon) \: c_2 \: & = \: 
1 \: - \: \beta, \notag \\
c_1 \: (1 \: + \: R)^{- \epsilon} \: + \:
c_2 \: (1 \: + \: R)^{1 \: + \: \epsilon} \: & = \: C \: (1 \: + \: \beta) \: 
\delta_0, \notag \\
- \: \epsilon \: c_1 \: (1 \: + \: R)^{- \epsilon - 1} \: + \:
(1 \: + \: \epsilon) \: 
c_2 \: (1 \: + \: R)^{\epsilon} \: & = \: 1 \: + \: \beta. 
\end{align}
(The third equation in (\ref{eq4.3}) says that the sizes of the circles
$\{R\} \times S^1$ and $C \cdot \partial_0(N_1)$ are the same,
while the fourth equation in (\ref{eq4.3}) says that the cone angles along the
circles are the same.)
The solution to the first two equations in (\ref{eq4.3}) is
\begin{align} \label{eq4.4}
c_1 \: & = \: \frac{1-\beta}{1+2\epsilon} \: \left( (1 \: + \: \epsilon) \:
\delta_\infty \: - \: 1 \right), \notag \\
c_2 \: & = \: \frac{1-\beta}{1+2\epsilon} \: \left( \epsilon \:
\delta_\infty \: + \: 1 \right).
\end{align}
For small $\epsilon$ and large $R$, the dominant term on the left-hand-side of
the last equation in (\ref{eq4.3}) is
$(1 \: + \: \epsilon) \: 
c_2 \: (1 \: + \: R)^{\epsilon}$. Hence for small $\epsilon$, there is a
solution for $R$ with the asymptotics
\begin{equation} \label{eq4.5}
R \sim \left(\frac{1+\beta}{1-\beta} \right)^{\frac{1}{\epsilon}}.
\end{equation}
Substituting into the third equation of (\ref{eq4.3}) gives
\begin{equation} \label{eq4.6}
C \sim \delta_0^{-1} \:
\left(\frac{1+\beta}{1-\beta} \right)^{\frac{1}{\epsilon}}.
\end{equation}

Put $P \: = \: N_1 \cup_{ S^1 } N_2$,
where the gluing identifies $\partial_\infty N_1$ with $\{0\} \times S^1
\subset N_2$.
Then $P$ has a $C^1$-smooth Riemannian metric which is flat on $N_1$ and has
curvature 
$-\: \frac{f^{\prime \prime}}{f} \: = \:
- \: \frac{\epsilon \: (1 \:  + \: \epsilon)}{(1 + r)^2}$ on $N_2$.
By smoothing the metric on $P$ and slightly moving the boundary curve
between $N_1$ and $N_2$ into $N_2$, we can construct a Riemannian metric
on $P$ which is flat on $N_1$, which satisfies
$|K| \: \le \: \frac{2 \: \epsilon \: (1 \:  + \: \epsilon)}{(1 + r)^2}$ on
$N_2$ and for which $P$ glues isometrically onto $C \cdot P$ by identifying
$\{R\} \times S^1 \subset P$ with $C \cdot \partial_0 N_1 \subset C \cdot P$. 
Let $D^2$ be a $2$-disk which caps $P$ at $\partial_0 N_1$.
Put 
\begin{equation} \label{eq4.7}
M \: = \: D^2 \: \cup_{S^1} P \: \cup_{S^1}  (C \cdot P) \cup_{S^1}
\ldots \cup_{S^1} (C^k \cdot P) \cup_{S^1} \ldots,
\end{equation}
with basepoint $\star \in D^2$.
There is an obvious Riemannian metric on $M - D^2$, which we extend over $M$.
We claim that this Riemannian metric satisfies the conditions of the theorem.
First, $M$ has infinite topological type. By the self-similar nature of the
Riemannian metric, equations (\ref{eq1.7}) and (\ref{eq1.8}) are satisfied for 
some $c, c^\prime_1, c^\prime_2 \: > \: 0$. In order to check (\ref{eq1.9})
on $(C^k \cdot P) \subset M$, 
we can use the scale invariance to instead check it on
the subset $P$ of
\begin{equation} \label{eq4.8}
C^{-k} \cdot M \: = \: (C^{-k} \cdot D^2) \: \cup_{S^1} (C^{-k} \cdot P) \: 
\cup_{S^1}  (C^{-k+1} \cdot P) \cup_{S^1}
\ldots \cup_{S^1} P \cup_{S^1} \ldots
\end{equation}
As the metric is flat on $N_1 \subset P$, it is enough to just consider 
a point $m \in N_2$, say with coordinates $(r, \theta) \in [0, R] \times S^1$.
Put 
\begin{equation} \label{eq4.9}
a_1 \: = \: \max_{z_1 \in \partial_0 N_1, z_2 \in \partial_\infty N_1}
d(z_1, z_2)
\end{equation}
and
\begin{equation} \label{eq4.10}
a_2 \: = \: \max_{z \in \partial D^2} d(\star, z).
\end{equation}
Then we can construct a path from $\star$ to $m$ with
length at most
\begin{equation} \label{eq4.11}
C^{-k} \: a_2 \: + C^{-k} (a_1 \: + \: R) \: + \: \ldots \: + \:
C^{-1} (a_1 \: + \: R) \: +
\: a_1 \: + \: r.
\end{equation}
Thus
\begin{equation} \label{eq4.12}
d(m, \star) \: \le \: a_2 \: + \: \frac{a_1 \: + \: R}{C \: - \: 1} \: + \:
a_1 \: + \: r
\: \le \: r \: + \: \const,
\end{equation}
where $\const$ is independent of $\epsilon$. It follows that
\begin{equation} \label{eq4.13}
|K(m)| \cdot d(m, \star)^2 \: \le \: 
2 \: \epsilon \: (1 + \epsilon) \: 
\left( \frac{r \: + \: \const}{r \: + \: 1} \right)^2,
\end{equation}
which proves the theorem. \\ \\
\noindent
{\bf Remark :} It should be fairly clear that by using building blocks
consisting of appropriate (rescaled) flat metrics on
$T^2  -  (D^2 \cup D^2)$, $S^2  - 
(D^2 \cup D^2 \cup D^2)$ and $\R P^2  -  (D^2 \cup D^2)$,
along with the classification of surfaces in \cite{Richards (1963)}, we
can construct a complete Riemannian metric on {\it any} connected surface
so as to satisfy (\ref{eq1.7}), (\ref{eq1.8}) and (\ref{eq1.9}) for any
$\epsilon \: > \: 0$ and for some $c$, $c^\prime_1$ and $c^\prime_2$.

\section{Proof of Theorem 3} \label{Proof of Theorem 3}

We follow the method of proof of Theorem \ref{thm1}, which is a proof
by contradiction. Hence we obtain a
pointed length space $(X_\omega, \star_\omega)$ along with a flat 
$n$-dimensional Riemannian
metric on $X_\omega - \star_\omega$. By using appropriate rescalings in
the construction of $X_\omega$, we obtain the analog of equations
(\ref{eq1.10}) and (\ref{eq1.11}) for 
$X_\omega$, but with $C^\prime \: = \: 0$. That is,
the distance 
function $d_{X_\omega}(\cdot, \star_\omega)$ is convex on $X_\omega$ and 
for any two normalized minimizing geodesics 
$\gamma_1, \gamma_2 \: : \:
[0,b] \rightarrow X_\omega$ with 
$\gamma_1(0) \: = \: \gamma_2(0) \: = \: \star_\omega$ and any $t \in [0,1]$,
\begin{equation} \label{eq5.1}
d_{X_\omega}(\gamma_1(t b), \gamma_2(tb)) \: \le \: t \: 
d_{X_\omega}(\gamma_1(b), \gamma_2(b))
\end{equation}

Let $c_\omega \in X_\omega$ and ${\cal C} \subset X_\omega - \star_\omega$
be as in the proof of Theorem \ref{thm1}. Let $\widetilde{\cal C}$ denote the
universal cover of ${\cal C}$, defined with the basepoint $c_\omega$, with
projection $\pi \: : \: \widetilde{\cal C} \rightarrow {\cal C}$.
As ${\cal C}$ is flat,
there is a developing map $D \: : \: \widetilde{\cal C} \rightarrow
\R^n$ and a homomorphism $\pi_1({\cal C}, c_\omega) \rightarrow
\Isom(\R^n)$ with respect to which $D$ is equivariant. 

From the convexity of $d(\cdot, \star_\omega)$, for any $r \: > \: 0$ the
ball $B_r(\star_\omega)$ is geodesically convex in ${\cal C}$. Then
$S_r(\star_\omega)$ is locally convex in the sense that for each 
${x}_\omega \in S_r(\star_\omega)$, there is a 
neighborhood of ${x}_\omega$ in $S_r(\star_\omega)$
which is contained in the boundary of a convex set. Given
$\widetilde{x}_\omega \in \pi^{-1}(S_r(\star_\omega))$,
using a local
isometry between a neighborhood of $\widetilde{x_\omega}$ and a neighborhood of
$\pi(\widetilde{x_\omega})$, it follows that there is a 
neighborhood of $\widetilde{x}_\omega$ in $\pi^{-1}(S_r(\star_\omega))$
which is contained in the boundary of a convex set.
That is, $\pi^{-1}(S_r(\star_\omega))$ is locally convex. From
\cite{Jonker-Norman (1973)}, for each $r \: > \: 0$, \\
1. $\pi^{-1}(S_r(\star_\omega))$ is 
embedded by $D$ as the boundary of a convex subset of
$\R^n$, or\\
2. $\pi^{-1}(S_r(\star_\omega))$ is 
isometric to $S^1 \times \R^{n-2}$ and $D$ is the product $\alpha \times
\Id_{\R^{n-2}}$ of an immersed convex curve $\alpha \: : \: S^1 \rightarrow
\R^2$ with the identity map on $\R^{n-2}$, or \\
3. $\pi^{-1}(S_r(\star_\omega))$ is 
isometric to $\R \times \R^{n-2}$ and $D$ is the product $\alpha \times
\Id_{\R^{n-2}}$ of an immersed convex curve $\alpha \: : \: \R \rightarrow
\R^2$ with the identity map on $\R^{n-2}$.

Suppose first that for each $r \: > \: 0$, $D$ embeds
$\pi^{-1}(S_r(\star_\omega))$ into $\R^n$ as the boundary of a convex subset.
Then $D$ is an embedding of $\widetilde{\cal C}$ into
$\R^n$. Identifying $\widetilde{\cal C}$ with its image under $D$, 
convexity implies that 
$\widetilde{\cal C}$ is the complement of a closed
convex subset $Z \subset \R^n$. Letting $\R^n/Z$ denote the collapsing of
$Z$ to a point, there is a continuous map $\R^n \rightarrow \R^n/Z 
\rightarrow \overline{\cal C}$ which sends $Z$ to $\star_\omega$.
Now $Z$ is invariant under
the isometric action of $\pi_1({\cal C}, c_\omega)$ on $\R^n$. 
Given $x_\omega \in {\cal C}_\omega$ and a lift $\widetilde{x}_\omega \in
\pi^{-1}(x_\omega)$, the convexity of $Z$ implies that there is a unique
minimizing geodesic from $x_\omega$ to $\star_\omega$, which coincides
with the projection of the minimizing segment from $\widetilde{x}_\omega$ to
$Z$. 

Suppose that $Z$ contains more than one point.  Then we can find two
distinct points
$\{z_i\}_{i = 1,2}$ in $\partial Z$ and support planes $H_i$
containing $z_i$ so that the normalized rays
$\{\widetilde{\gamma}_i\}_{i = 1,2}$ from $z_i$ orthogonal to $H_i$, which
point away from $Z$, have the property that $\widetilde{\gamma}_1$
eventually lies on the same side of $H_2$ as $\widetilde{\gamma}_2$, 
and  $\widetilde{\gamma}_2$
eventually lies on the same side of $H_1$ as $\widetilde{\gamma}_1$.
Put $\gamma_i \: = \: \pi \circ 
\widetilde{\gamma}_i$. For $t$ sufficiently small, we will
have $d_{X_\omega}(\gamma_1(t), \gamma_2(t)) \: = \: 2t$, 
as the shortest way to
get from $\widetilde{\gamma}_1(t)$ to $\widetilde{\gamma}_2(t)$
in $\R^n/Z$ will
be to follow $\widetilde{\gamma}_1$ from $\widetilde{\gamma}_1(t)$ to
$z_1$ and then follow $\widetilde{\gamma}_2$ from $z_2$ to 
$\widetilde{\gamma}_2(t)$. (Note that $Z$ gets collapsed to $\star_\omega$.)
Then from (\ref{eq5.1}), it follows that 
$d_{X_\omega}(\gamma_1(t), \gamma_2(t)) \: = \: 2t$ for all $t \: > \: 0$. Thus
$d(\widetilde{\gamma}_1(t), \widetilde{\gamma}_2(t)) \: = \: 2t$ for all 
$t \: > \: 0$, where the distance is measured in the length metric on 
$\R^n/Z$, which is a contradiction to the construction of 
$\widetilde{\gamma}_1$ and $\widetilde{\gamma}_2$.

Thus $Z$ is a point, which we can assume 
without loss of generality to be the origin in $\R^n$. 
Then 
$\pi_1({\cal C}, c_\omega)$ acts on $\R^n - \{0\}$ by elements of $O(n)$ and
${\cal C}$ is a cone over a spherical space form.  The rest of the
proof proceeds as in the proof of Theorem \ref{thm1}.

Now suppose that for some $r_0 \: > \: 0$, $D$ immerses
$\pi^{-1}(S_{r_0}(\star_\omega))$ as
$\alpha_{r_0} \times \Id_{\R^{n-2}}$ , where 
$\alpha_{r_0}$ is an immersed convex curve 
$\alpha_{r_0} \: : \: S^1 \rightarrow \R^2$. Then for all $r \: > \: 0$,
$D$ immerses
$\pi^{-1}(S_{r}(\star_\omega))$ as
$\alpha_{r} \times \Id_{\R^{n-2}}$ , where 
$\alpha_{r}$ is an immersed convex curve 
$\alpha_r \: : \: S^1 \rightarrow \R^2$ which is the curve of distance
$r \: - \: r_0$ from $\alpha_{r_0}$. (Recall that $D$ is a local
isometry.)  It follows that 
$\widetilde{\cal C} \: = \: (0, \infty) \times S^1 \times 
\R^{n-2}$, which contradicts the fact that $\widetilde{\cal C}$ is
simply-connected.

Finally, suppose that for some $r_0 \: > \: 0$, $D$ immerses
$\pi^{-1}(S_{r_0}(\star_\omega))$ as
$\alpha_{r_0} \times \Id_{\R^{n-2}}$ , where 
$\alpha_{r_0}$ is an immersed convex curve 
$\alpha_{r_0} \: : \: \R \rightarrow \R^2$. Then for all $r \: > \: 0$,
$D$ immerses
$\pi^{-1}(S_{r}(\star_\omega))$ as
$\alpha_{r} \times \Id_{\R^{n-2}}$ , where 
$\alpha_{r}$ is an immersed convex curve 
$\alpha_r \: : \: \R \rightarrow \R^2$ which is the curve of distance
$r \: - \: r_0$ from $\alpha_{r_0}$.
In particular, $\widetilde{\cal C}$ splits isometrically as a product
${\cal A} \times \R^{n-2}$, where
${\cal A}$ is diffeomorphic to $(0,\infty) \times \R$, 
with $\widetilde{\cal C}$ having the flat metric which pulls back from $D$. Put
$\overline{\cal A} \: = \: ([0, \infty) \times \R)/(\{0\} \times \R)$,
the union of ${\cal A}$ with a point.  Similarly, put
$\overline{\widetilde{\cal C}} 
\: = \: ([0, \infty) \times \R \times \R^{n-2})/(\{0\} \times
\R \times \R^{n-2})$, 
the union of $\widetilde{\cal C}$ with a point. There is a continuous
map from $\overline{\widetilde{\cal C}}$ to $\overline{\cal C}$ which
restricts to the covering map on $\widetilde{\cal C}$, and an obvious
embedding
$\overline{\cal A} \rightarrow \overline{\widetilde{\cal C}}$.
Let
$\widetilde{\gamma} : [0, \infty) \rightarrow \overline{\cal A}$ be a
normalized ray. Choose distinct points
$b_1, b_2 \in \R^{n-2}$. Then $(t \in \R^+) \rightarrow 
\widetilde{\gamma}(t) \times \{b_1\}$ and
$(t \in \R^+) \rightarrow 
\widetilde{\gamma}(t) \times \{b_2\}$ extend to rays
$r_i \: : \: [0, \infty) \rightarrow \overline{\widetilde{\cal C}}$,
with $r_i(0)$ being the basepoint. As before, we have
$d(r_1(t), r_2(t)) \: = \: 2t$ for $t$ small, where $d$ is the
length metric on $\overline{\widetilde{C}}$.  Then (\ref{eq5.1}) implies
that $d(r_1(t), r_2(t)) \: = \: 2t$ for all $t$, which is a 
contradiction. \\ \\
\noindent
{\bf Remark : } To see where the hypotheses of Theorem \ref{thm3} 
enter into the proof,
note that the method of proof is to show that
$\overline{\cal C}$ is a cone over a spherical space form. 
If $n \: = \: 2$ then $\overline{\cal C}$ could {\it a priori} be a cone
over $\R$, as in the example of Section \ref{Ultralimits}. To see where the
assumption of large-scale pointed-convexity enters, let $M_i$ be the
effect of attaching a wormhole between two points of distance
$2i$ in $\R^n$. More precisely, give $[-1/2, 1/2] \times S^{n-1}$ a metric
whose restrictions to $[-1/2, -1/4] \times S^{n-1}$ and 
$[1/4, 1/2] \times S^{n-1}$ are isometric to $\overline{B_{1/2}(0)} -
{B_{1/4}(0)} \subset \R^n$. Put $M_i \: = \: 
(\R^n - B_{1/2}(p_1) - B_{1/2}(p_2)) \cup_{S^{n-1} \cup S^{n-1}} 
[-1/2, 1/2] \times S^{n-1}$, where $p_1, p_2 \in \R^n$ have distance
$2i$. It is flat outside of a compact set. Put the basepoint of $M_i$ somewhere
on the component $[-1/2, 1/2] \times S^{n-1}$.
Then the limit space $X_\omega \: = \: 
\lim_\omega \frac{1}{i} \cdot M_i$ is the result of identifying two points in 
$\R^n$ of distance $2$, with its basepoint at the identification point.
Clearly $X_\omega$ is not a cone. 
Without the assumption of large-scale pointed-convexity, 
or some such assumption, it could {\it a priori} arise in a
rescaling limit as in the proof of Theorem \ref{thm3}. One can find similar
examples with $X_\omega \: = \: \R^n/K$, where $K$ is any closed subset of
$\R^n$. The large-scale pointed-convexity assumption is used to show
first that $K$ is convex and then to show that $K$ is a point.

\end{document}